\documentclass[12pt]{article}

\usepackage{amsmath,amsthm,amsfonts}

\parskip=\medskipamount
\def\onehalf{{\textstyle\frac12}}

\def\cinfty#1{C^{\scriptscriptstyle\infty}(#1)}
\def\vectorfields#1{{\cal X}(#1)}

\def\ov#1{\overline{#1}}
\def\del{{\nabla}}
\def\lie#1{{\cal L}_{#1}}
\def\fpd#1#2{\frac{\partial #1}{\partial #2}}
\def\R{{\rm I\kern-.20em R}}

\def\proof{{\sc Proof.}}

\newtheorem{prop}{\bf Proposition}
\newtheorem{lemma}{\bf Lemma}
\newtheorem{thm}{\bf Theorem}

\def\T{{\bf T}}
\def\DV#1{\V{\rm D}_{#1}}
\def\DH#1{\H{\rm D}_{#1}}
\def\dv#1{d^{\scriptscriptstyle V}_{#1}}
\def\dh#1{d^{\scriptscriptstyle H}_{#1}}
\def\H#1{{#1}^{\scriptscriptstyle H}}
\def\V#1{{#1}^{\scriptscriptstyle V}}
\def\K#1{{#1}^{\scriptscriptstyle K}}
\def\subH#1{{#1}_{\scriptscriptstyle H}}
\def\subV#1{{#1}_{\scriptscriptstyle V}}

\def\hook{\vrule height 0pt depth 0.4pt width 3pt \vrule height 7pt depth 0.4pt
\kern 3pt}

\begin{document}

\title{A class of Poisson-Nijenhuis structures on a tangent bundle}

\author{W.\ Sarlet and F.\ Vermeire\footnote{Research Assistant of the Fund for Scientific Research-Flanders, Belgium (FWO-Vlaanderen)}\\
{\small Department of Mathematical Physics and Astronomy }\\
{\small Ghent University, Krijgslaan 281, B-9000 Ghent, Belgium}}
\date{}

\maketitle

\begin{quote}
{\bf Abstract.} {\small Equipping the tangent bundle $TQ$ of a
manifold with a symplectic form coming from a regular Lagrangian
$L$, we explore how to obtain a Poisson-Nijenhuis structure from a
given type $(1,1)$ tensor field $J$ on $Q$. It is argued that the
complete lift $J^c$ of $J$ is not the natural candidate for a
Nijenhuis tensor on $TQ$, but plays a crucial role in the
construction of a different tensor $R$, which appears to be the
pullback under the Legendre transform of the lift of $J$ to
$T^*Q$. We show how this tangent bundle view brings new insights
and is capable also of producing all important results which are
known from previous studies on the cotangent bundle, in the case
that $Q$ is equipped with a Riemannian metric. The present
approach further paves the way for future generalizations.}
\end{quote}

\section{Introduction}
There is a well established theory of bi-Hamiltonian systems on the cotangent
bundle $T^*Q$ of a manifold $Q$ or, more generally, on a symplectic or Poisson
manifold. In particular, there is a link between separability of the
Hamilton-Jacobi equation and certain classes of bi-Hamiltonian or
quasi-bi-Hamiltonian systems, in which more specifically Poisson-Nijenhuis
structures \cite{KosMa} play a prominent role. The immediate source of
inspiration for the present work is a series of recent papers in this general
field in which relations have been explored between such things as:
bi-differential calculi, complete integrability, St\"ackel systems, compatible
Poisson structures on an extended space, Gelfand-Zakharevich systems,
so-called special conformal Killing tensors, and so on. See
\cite{Blaszak,Blaszak2,Crampin,CraSar1,CraSar2,CST,Falqui,FP,%
GelZak,IMM,Lundmark,MaBla,MorTon} for a non-exhaustive list of
recent contributions. In the more direct applications of these
theoretical developments, the dominant geometrical space is a
cotangent bundle, which of course comes equipped with its
canonical symplectic or Poisson structure, and on which a
compatible Poisson structure is obtained via the lift $\tilde J$
of a type $(1,1)$ tensor field $J$ on the base manifold, thus
giving birth to a Poisson-Nijenhuis structure. The link with a
class of St\"ackel systems requires the availability of a metric
$g$ on $Q$, with respect to which $J$ has the property of being a
special conformal Killing tensor (a concept which is intimately
linked with what is called Benenti tensor \cite{Blaszak2,BolMat},
after the extensive work of Benenti on Hamilton-Jacobi
separability (see e.g.\ \cite{Ben,Ben2})).

The point we want to emphasize now, however, is that some of these
applications clearly come from meaningful questions about dynamical systems
living on a tangent bundle $TQ$. This is for example the case with the theory
of cofactor pair systems, as developed originally by Lundmark on a Euclidean
space \cite{Lundmark,Lundmark2} and generalised to (pseudo-)Riemannian
manifolds and in more geometrical terms by Crampin and Sarlet \cite{CraSar1}.
The physical background for these systems is a kinetic energy type Lagrangian
on $TQ$ for which (at least in one possible interpretation) admissible
non-conservative forces are being sought, in the sense that the resulting
Newtonian system admits two quadratic first integrals, which in turn can
generate a whole family of integrals in involution. The cofactors of the
Killing tensors coming from these integrals then determine the special
conformal Killing tensors which give rise to Poisson-Nijenhuis structures; we
are then looking at examples of so called bi-quasi-Hamiltonian systems
\cite{CraSar2}.

A type $(1,1)$ tensor field $J$ on $Q$ also has a natural lift to $TQ$; it is
usually called the {\sl complete lift\/} and we will denote it by $J^c$. The
difference with $T^*Q$ of course is that $TQ$ does not carry a canonical
Poisson structure. However, a symplectic form is available (and can be
constructed by pure tangent bundle techniques) as soon as a regular Lagrangian
is given, which could for example be the kinetic energy Lagrangian coming from
a metric on $Q$. It does, therefore, perfectly make sense to explore the
possibility of obtaining Poisson-Nijenhuis structures on $TQ$ by natural
tangent bundle constructions. One of the goals we have in mind for the future
is to arrive at a generalization of the theory of special conformal Killing
tensors from Riemannian to Finsler spaces. The primary objective of this
paper, however, is to set the ground for future developments by trying to
understand in detail how the results one is by now familiar with in a
cotangent bundle environment, can be obtained in a natural way by pure tangent
bundle techniques. We shall see that this different way of approaching the
subject offers new insights anyway. In fact, some of the preliminary
considerations lead to results which are valid for arbitrary Lagrangians, not
just kinetic-energy type ones. In this respect, we are to some extent joining
the interest in Poisson structures on a tangent bundle which is present also
in recent work by Vaisman \cite{MV,Vaisman1,Vaisman2}.

In Section 2, starting from a given Lagrangian $L$ on $TQ$ and a type $(1,1)$
tensor field $J$ on $Q$, we show that $J^cS$ determines an alternative almost
tangent structure and use this to construct another type $(1,1)$ tensor $R$ on
$TQ$. The full characterization of $R$ is developed in Section~3, making use
of the connection provided by the second-order equation field $\Gamma$, coming
from the Lagrange equations. The eigenspace structure of $R$ is discussed in
Appendix~B, whereas further properties of general interest are derived in
Section~4. In Section~5, we specialize to the particular case of a Riemannian
manifold and the associated kinetic energy Lagrangian and explore how various
known properties make their appearance within such a tangent bundle approach.
Some unexpected new features come forward which are further discussed in
Appendix~A. The road map to future developments is sketched in the final
section.

\section{Generalities}

Suppose we are given a regular Lagrangian $L$ on $TQ$ and a type $(1,1)$
tensor field $J$ on $Q$. $L$ comes with its associated symplectic form
$\omega_L$, and $J$ determines a tensor field $J^c$, its {\sl complete
lift\/}, on $TQ$. One may wonder whether these data can give rise, under some
circumstances, to a compatible Poisson structure. For that, $J^c$ should be
symmetric with respect to $\omega_L$, the so-called Magri-Morosi concomitant
must vanish (see \cite{NunezMarle,MaMor}), and also the Nijenhuis torsion of
$J^c$ must be zero (which is equivalent to the torsion of $J$ being zero). It
seems to us, however, that this is not the most interesting path to pursue.
Experience in a variety of applications (see for example
\cite{Clift,MFLMR,SaCra}) has shown that interesting type $(1,1)$ tensor
fields $R$ on a symplectic manifold $(M,\omega)$ arise from the construction
of a second 2-form $\omega_1$ and the determining formula:
\begin{equation}
i_{R(\xi)}\omega = i_\xi\omega_1, \qquad \forall
\xi\in\vectorfields{TQ}. \label{generalR}
\end{equation}
A direct advantage of such tensor fields is for example that they
automatically have the required symmetry property
\begin{equation}
\omega(R\xi,\eta) = \omega(\xi,R\eta). \label{Rsymmetry}
\end{equation}
Furthermore, vanishing of the Magri-Morosi concomitant then is
equivalent to $d\omega_1=0$ \cite{CST}, after which we are left
with the condition $N_R=0$. Although the canonical lift
$\tilde{J}$ of $J$ to $T^*Q$ can be defined (see \cite{Clift}) via
a relation like (\ref{generalR}), this does not seem to be the
case for $J^c$ on $TQ$. We therefore start our investigation with
an exploration of possible natural ways of constructing a second
closed 2-form $\omega_1$ from the given data on $TQ$.

For completeness, we list a number of useful properties of $J^c$,
starting with defining relations with respect to the action on
complete and vertical lifts of vector fields on $Q$: for $X,Y\in
\vectorfields{Q}$,
\begin{equation}
J^c(X^c)=  (JX)^c, \qquad  J^c(\V{X})= \V{(JX)}. \label{defJc}
\end{equation}
Using the bracket relations
\[
{}[\V{X},\V{Y}]=0, \quad [\V{X},Y^c]=\V{[X,Y]}, \quad
[X^c,Y^c]=[X,Y]^c,
\]
it easily follows that
\begin{equation}
\lie{X^c}J^c= (\lie{X}J)^c, \qquad \lie{\V{X}}J^c=\V{(\lie{X}J)}.
\label{lieJc}
\end{equation}
Also, for the Nijenhuis torsion we have,
\begin{eqnarray*}
N_{J^c}(\V{X},\V{Y})&=& 0,\\
N_{J^c}(X^c,Y^c)&=& (N_J(X,Y))^c, \\
N_{J^c}(\V{X},Y^c) &=& N_{J^c}(X^c,\V{Y}) \;=\; \V{(N_J(X,Y))},
\end{eqnarray*}
from which it follows that $N_{J^c}=0\, \Longleftrightarrow \,
N_J=0$.

It is imperative to relate $J^c$ to the canonical type $(1,1)$ tensor field
$S$ on $TQ$ (the so-called vertical endomorphism), which satisfies
$S(X^c)=\V{X},\ S(\V{X})=0$. It is easy to see that $J^c$ commutes with $S$,
in the sense of endomorphisms on $\vectorfields{TQ}$, but also in the sense of
the Nijenhuis bracket:
\begin{equation}
[J^c,S] = 0. \label{JcS}
\end{equation}
It follows that also the corresponding degree 1 derivations
commute, meaning that:
\begin{equation}
d_S\,d_{J^c} = - d_{J^c}\,d_S. \label{dSdJc}
\end{equation}
Finally, it is easy to verify that for any $J$,
\begin{equation}
N_{J^cS}=0 . \label{NJcS}
\end{equation}
In fact, we can make the following more complete statement in that
respect, which is trivial to prove, and essentially says that
$J^cS$ very much has the same properties as $S$.

\begin{lemma} We have $(J^cS)^2=0$ and $N_{J^cS}=0$. Furthermore, if $J$ is
non-singular, $J^cS$ determines an integrable almost tangent
structure.
\end{lemma}

Now recall the role which the vertical endomorphism plays in the
definition of the Poincar\'e-Cartan 2-form $\omega_L$: we have,
\[
\omega_L = d(S(dL)) = dd_SL.
\]
[We make no notational distinction between the action of a type
$(1,1)$ tensor field on vector fields and its dual action on
1-forms, but one should keep in mind that the order of composition
of such action changes in passing from one interpretation to the
other.] It then looks perfectly natural, given $J$ and the sort of
alternative integrable almost tangent structure which it creates,
to consider the closed 2-from $\omega_1$, defined (with various
ways of writing the same expression) by,
\begin{equation}
\omega_1= d(SJ^c(dL))=d(S(d_{J^c}L)) = d(J^c(d_S L))= dd_{J^cS}L.
\label{omega1}
\end{equation}
And so, with $L$ and $J$ as data, the type $(1,1)$ tensor field
$R$ which will carry our attention is defined by
\begin{equation}
i_{R(\xi)}dd_SL = i_\xi dd_{J^cS}L. \label{ourR}
\end{equation}
We know that it will define a Poisson-Nijenhuis structure if and
only if $N_R=0$.

The first objective now must be to obtain a reasonably practical
description of $R$, for example by recognising its action on
complete and vertical lifts. In fact, we believe that it is better
for general purposes, to make use of horizontal and vertical
lifts, rather than complete and vertical lifts. For that, of
course, one needs a connection, but there is one available, namely
the non-linear connection associated to the Euler-Lagrange
equations of $L$ (being second-order differential equations on
$TQ$).

\section{Making use of a connection}

As is well known, every second-order equation field
\begin{equation}
\Gamma = u^i\fpd{}{q^i} + f^i(q,u)\fpd{}{u^i}, \label{sode}
\end{equation}
in particular the one coming from the regular Lagrangian $L$, defines a
horizontal distribution, with connection coefficients
\begin{equation}
\Gamma^i_j = - \onehalf \fpd{f^i}{u^j}. \label{conn}
\end{equation}
As a result, every vector field $\xi$ on $TQ$ has a unique decomposition of
the form $\xi=\H{X}+\V{Y}$, where $X,Y$ are vector fields along the tangent
bundle projection $\tau:TQ \rightarrow Q$. The $\cinfty{TQ}$ module of such
fields is denoted by $\vectorfields{\tau}$. An extensive calculus along the
projection $\tau$ was developed in \cite{MaCaSaI,MaCaSaII}. We recall here
some basic features of this calculus, which will be needed in what follows.

Interesting derivations and tensorial objects along $\tau$ are discovered by
looking at the decomposition of Lie brackets of vector fields on $TQ$. We
have, for example, that
\begin{eqnarray}
[\V{X},\V{Y}] &=& \V{(\subV{[X,Y]})}, \label{[vv]} \\
{}[\H{X},\V{Y}] &=& \V{(\DH{X}Y)} - \H{(\DV{Y}X)}, \label{[hv]} \\
{}[\H{X},\H{Y}] &=& \H{(\subH{[X,Y]})} + \V{{\mathcal R}(X,Y)}. \label{[hh]}
\end{eqnarray}
Here, ${\mathcal R}$ is the curvature of the non-linear connection, a vector
valued 2-form along $\tau$, $\DH{X}$ and $\DV{X}$ are the horizontal and
vertical covariant derivative operators, which act on functions
$F\in\cinfty{TQ}$ as
\begin{equation}
\DH{X}F =\H{X}(F), \qquad \DV{X}F = \V{X}(F). \label{DHandDV}
\end{equation}
The horizontal and vertical brackets of elements of
$\vectorfields{\tau}$ (in the present situation of a connection
which has no torsion) are given by
\begin{equation}
\subV{[X,Y]}= \DV{X}Y - \DV{Y}X, \quad \subH{[X,Y]}=\DH{X}Y-\DH{Y}X .
\label{vhbrackets}
\end{equation}
Other brackets of interest are
\begin{equation}
[\Gamma,\V{X}]= -\H{X} + \V{(\del X)} , \qquad [\Gamma,\H{X}]=  \H{(\del X)} +
\V{\Phi(X)}. \label{gammabrackets}
\end{equation}
Here $\Phi$, a type $(1,1)$ tensor along $\tau$, is called the
{\sl Jacobi endomorphism\/}, and $\del$ is the {\sl dynamical
covariant derivative\/}, which on functions acts like $\Gamma$.

There exist lift operations on many other tensorial objects. We mention two
more constructions of interest now. First there are the horizontal and
vertical lifts of a type $(1,1)$ tensor field $U$ along $\tau$, determined by
\begin{eqnarray}
\H{U}(\V{X})=\V{U(X)}, &\quad& \H{U}(\H{X})=\H{U(X)}, \label{uh} \\
\V{U}(\V{X})=0, &\quad& \V{U}(\H{X})=\V{U(X)}. \label{uv}
\end{eqnarray}
Next, if $g$ is a symmetric type $(0,2)$ tensor field along $\tau$, its
K\"ahler lift is a 2-form on $TQ$ determined by
\begin{eqnarray}
&& \K{g}(\H{X},\H{Y}) = \K{g}(\V{X},\V{Y}) =0, \label{gk1} \\
&& \K{g}(\V{X},\H{Y}) = - \K{g}(\H{X},\V{Y}) = g(X,Y). \label{gk2}
\end{eqnarray}
In fact, the Poincar\'e-Cartan form $\omega_L$ is precisely the K\"ahler lift
of the Hessian of $L$, defined intrinsically by
\begin{equation}
g = \DV{}\DV{}L, \label{dvdvL}
\end{equation}
where $\DV{}$ is the vertical covariant differential defined (on
any tensor $T$) by $\DV{}T(X,\ldots)=\DV{X}T(\ldots)$. Observe
that the complete lift $J^c$ of a $(1,1)$ tensor $J$ on $Q$, can
be written as
\begin{equation}
J^c = \H{J} + \V{(\del J)}. \label{Jc}
\end{equation}

Other, more specific properties of interest in the calculus along
$\tau$ will be recalled when appropriate, but we should at least
refer here also to the existence of a canonical vector field along
$\tau$, the {\sl total time derivative\/} operator
$\T=u^i\partial/\partial q^i$, whose vertical lift is the
Liouville vector field $\Delta$, whereas its horizontal lift is a
second-order equation field, which need not be the one we started
from; the two coincide for a quadratic spray.

We can now start the computation of the structure of the tensor field $R$
defined by (\ref{ourR}). To begin with, observe that from
(\ref{uh}),(\ref{uv}) and (\ref{Jc}), it easily follows that
\begin{eqnarray}
J^c(\V{X}) &=& \V{(JX)}, \label{Jcxv} \\
J^c(\H{X}) &=& \H{(JX)} + \V{\del J(X)}. \label{Jcxh}
\end{eqnarray}
Using the standard notation $\theta_L$ for $d_SL$, we have
\begin{equation}
\theta_L(\V{X})=0, \qquad \theta_L(\H{X})= dL(\V{X})=\V{X}(L)
=\DV{X}L=\dv{}\!L(X). \label{thetaL}
\end{equation}
Note that in the second of these relations, a computation on $TQ$
in the end is replaced by one involving a vector field and 1-form
along $\tau$. The point is that $\theta_L$, being semi-basic, can
be regarded as a 1-form along $\tau$ as well; the defining
relation then reads $\theta_L=\dv{}\!L$, where $\dv{}$ is the
vertical exterior derivative. The latter is completely determined
by the following action on functions $F$ and 1-forms $\alpha$
along $\tau$:
\begin{equation}
\dv{}\!F(X)= \DV{X}F, \qquad \dv{}\!\alpha(X,Y)= \DV{X}\alpha(Y) -
\DV{Y}\alpha(X). \label{dv}
\end{equation}
Similar relations hold for the horizontal exterior derivative
$\dh{}$. The action of these exterior derivatives extends to
vector-valued forms as well; it suffices for our purposes to know
that for $\dv{}$ or $\dh{}$ acting on a type $(1,1)$ tensor field
along $\tau$, the defining relation is formally the same as in the
second of equations (\ref{dv}).

It is easy to see that $\omega_1=d(J^c\theta_L)$ gives zero when evaluated on
two vertical vector fields. Next we have, passing as before from $\theta_L$,
regarded as 1-form on $TQ$, to its interpretation as 1-form along $\tau$,
\begin{eqnarray*}
\omega_1(\V{X},\H{Y}) &=& \lie{\V{X}}\Big(\theta_L(\H{(JY)})\Big)
- \theta_L\Big(J^c([\V{X},\H{Y}])\Big) \\
&=& \DV{X}\Big(\theta_L(JY)\Big) - \theta_L\Big(J(\DV{X}Y)\Big) =
\DV{X}\theta_L(JY),
\end{eqnarray*}
from which it follows in view of (\ref{dvdvL}) that
\begin{equation}
\omega_1(\V{X},\H{Y}) = g(X,JY). \label{omega1vh}
\end{equation}
Proceeding in the same way, we get
\begin{eqnarray*}
\omega_1(\H{X},\H{Y}) &=& \lie{\H{X}}\Big(\theta_L(\H{(JY)})\Big)
-
\lie{\H{Y}}\Big(\theta_L(\H{(JX)})\Big) - \theta_L\Big(J^c([\H{X},\H{Y}])\Big) \\
&=& \DH{X}\Big(\theta_L(JY)\Big) - \DH{Y}\Big(\theta_L(JX)\Big) -
\theta_L\Big(J(\subH{[X,Y]})\Big) \\
&=& \DH{X}(J\theta_L(Y)) - \DH{Y}(J\theta_L(X)) -
\theta_L\Big(J(\DH{X}Y - \DH{Y}X)\Big) \\
&=& \DH{X}(J\theta_L)(Y) - \DH{Y}(J\theta_L)(X).
\end{eqnarray*}
It follows that
\begin{equation}
\omega_1(\H{X},\H{Y}) = \dh{}\!(J\theta_L)(X,Y). \label{omega1hh}
\end{equation}

\begin{prop} The type $(1,1)$ tensor field $R$, defined by (\ref{ourR}) has
the following structure:
\begin{eqnarray}
R(\V{X}) &=& \V{(\ov{J}X)} \label{Rvx} \\
R(\H{X}) &=& \H{(JX)} + \V{(UX)}, \label{Rhx}
\end{eqnarray}
where $\ov{J}$ is the transpose of $J$ with respect to
$g=\DV{}\DV{}L$, i.e.\ $g(JX,Y)=g(X,\ov{J}Y)$, and $U$ is the
tensor field along $\tau$, determined by
\begin{equation}
g(UX,Y) =  \dh{}\!(J\theta_L)(X,Y). \label{generalU}
\end{equation}
\end{prop}
\proof\ \ It is sufficient to take horizontal and vertical lifts of basic
vector fields for finding the tensorial structure of $R$. We have
$\omega_L(R(\V{X}),\V{Y})=0$ and
\[
\omega_L(R(\V{X}),\H{Y}) =  \omega_1(\V{X},\H{Y})= g(JY,X) =
g(Y,\ov{J}X) = \K{g}(\V{(\ov{J}X)},\H{Y}),
\]
from which (\ref{Rvx}) follows. Likewise
\[
\omega_L(R(\H{X}),\V{Y}) = - g(JX,Y) = \K{g}(\H{(JX)},\V{Y}),
\]
from which it follows that $R(\H{X})=\H{(JX)}+\V{(UX)}$, for some
$U$. Subsequently, using (\ref{gk2}) and (\ref{omega1hh}),
\[
\omega_L(R(\H{X}),\H{Y}) = \omega_L(\V{(UX)},\H{Y}) = g(UX,Y) =
\dh{}\!(J\theta_L)(X,Y),
\]
which completes the proof. \qed

Note that it follows from the skew-symmetry of the right-hand side
in (\ref{generalU}), that $\ov{U}=-U$.

In Appendix~B, we investigate the eigenspace structure of $R$ and the explicit
construction of so-called Darboux-Nijenhuis coordinates.

\section{Further properties of the tensor field $R$}

\begin{prop} We have $R=J^c\  \Longleftrightarrow\  J=\ov{J}
\mbox{\ \ and\ \ } U=\del J$.
\end{prop}
\proof\ \ The result follows immediately from comparison of
(\ref{Rvx}-\ref{Rhx}) with (\ref{Jcxv}-\ref{Jcxh}). \qed

A natural question which arises is whether $R$, in general, could
commute with $S$, just as $J^c$ does, either in the algebraic
sense or with respect to the Nijenhuis bracket.

\begin{prop} $RS=SR\ \Longleftrightarrow\ J=\ov{J}$.
\end{prop}
\proof\ \ Using $S(\V{X})=0$ and $S(\H{X})=\V{X}$, the result
follows immediately from the characterization of $R$ in
Proposition~1. \qed

Recall that the Nijenhuis bracket is defined by
\begin{eqnarray*}
[R,S](\xi,\eta)&=& [R\xi,S\eta] + [S\xi,R\eta] +
(RS+SR)([\xi,\eta]) \\
&& \mbox{} - R([S\xi,\eta]+[\xi,S\eta]) -
S([R\xi,\eta]+[\xi,R\eta]).
\end{eqnarray*}

\begin{prop} Assuming $J=\ov{J}$, so that $RS=SR$, we have
$[R,S]=0\ \Longleftrightarrow\ \dv{}\!U=\dh{}\!J$.
\end{prop}
\proof\ \ That $[R,S]$ vanishes on two vertical lifts is trivial. Again, it
suffices for such calculations to consider lifts of basic vector fields
(vector fields on $Q$), rather than vector fields along $\tau$. Since $J$ is
basic as well, one then easily verifies, making use of the bracket relations
(\ref{[vv]}) and (\ref{[hv]}), that also $[R,S](\V{X},\H{Y})=0$. Next, we have
\begin{eqnarray*}
[R,S](\H{X},\H{Y}) &=& [\H{(JX)}+\V{(UX)},\V{Y}] + [\V{X},\H{(JY)}+\V{(UY)}]  \\
&& \mbox{} + 2RS\Big(\H{\subH{[X,Y]}}\Big) -
R\Big([\V{X},\H{Y}]+[\H{X},\V{Y}]\Big) \\
&& \mbox{} - S\Big([\H{(JX)}+\V{(UX)},\H{Y}]+[\H{X},\H{(JY)}+\V{(UY)}]\Big).
\end{eqnarray*}
Using the bracket relations (\ref{[vv]}-\ref{[hh]}), plus the fact that $J$,
and by assumption also $X$ en $Y$, are basic, this readily reduces to
\[
[R,S](\H{X},\H{Y})= \V{(\DV{X}U(Y)-\DV{Y}U(X))} + \V{(J(\DH{X}Y-\DH{Y}X))} +
\V{(\DH{Y}(JX)-\DH{X}(JY))},
\]
which in view of properties such as (\ref{dv}) for type $(1,1)$ tensors, can
be written as
\[
[R,S](\H{X},\H{Y}) = \V{(\dv{}\!U(X,Y) - \dh{}\!J(X,Y))}.
\]
The result now follows. \qed

Another obvious question is under which circumstances $R$ is truly a recursion
operator for symmetries of $\Gamma$. For that, we compute $\lie{\Gamma}R$.
Taking the Lie derivative with respect to $\Gamma$ of the defining relations
(\ref{Rvx}-\ref{Rhx}) and making use also of the properties
(\ref{gammabrackets}), it is fairly straightforward to verify that
\begin{eqnarray}
\lie{\Gamma}R(\V{X}) &=& \H{(J - \ov{J})(X)} + \V{(U+\del\ov{J})(X)}, \label{lgamR1}\\
\lie{\Gamma}R(\H{X}) &=& \H{(\del J-U)(X)} + \V{(\del U + \Phi J-\ov{J}
\Phi)(X)}. \label{lgamR2}
\end{eqnarray}
The following, therefore, is an interesting immediate result.

\begin{prop}\ \ $\lie{\Gamma}R=0\ \Longleftrightarrow\ J=\ov{J},\ \ U=\del
J=0,\ \ \Phi J = J \Phi$. \qed
\end{prop}

Next, we address the question of recognizing the conditions under which $R$
has zero Nijenhuis torsion. A direct computation of $N_R$, through its action
on horizontal and vertical lifts is extremely tedious and therefore not worth
the effort, since we can actually rely on existing results concerning the
complete lift $\tilde{J}$ of $J$ to the cotangent bundle $T^*Q$. Let $Leg:
TQ\rightarrow T^*Q$ denote the Legendre transform coming from the regular
Lagrangian $L$.

\begin{prop}\ \ The tensor field $R$ defined by (\ref{ourR}) is directly
related to the complete lift $\tilde{J}$ on $T^*Q$: as a matter of
fact we have $R=Leg^*\tilde{J}$.
\end{prop}
\proof\ \  The defining relation of $\tilde{J}$, as described in \cite{Clift},
reads (where $\xi$ this time denotes an arbitrary vector field on $T^*Q$),
\begin{equation}
i_{\tilde{J}\xi}d\theta = i_\xi \lie{J^v}d\theta, \label{Jtilde}
\end{equation}
and clearly has the same structure as (\ref{ourR}). Here, $J^v$ is the
vertical lift of $J$, which is a vector field on $T^*Q$ (cf.\ \cite{Yano}) and
$\theta$ of course is the canonical 1-form on $T^*Q$. It is easy to see from
the coordinate expression that $i_{J^v}d\theta=\tilde{J}\theta$, so that
(\ref{Jtilde}) implies $i_{\tilde{J}\xi}d\theta = i_\xi d\tilde{J}\theta$. But
it is equally trivial to verify in coordinates that $Leg^*\tilde{J}\theta=
J^c\theta_L=d_{J^cS}L$. The result then immediately follows from taking the
pullback under $Leg$ of the new representation of (\ref{Jtilde}). \qed

As an immediate consequence of the fact that $N_{\tilde{J}}=0 \Leftrightarrow
N_J=0$, which incidentally requires also a fairly tedious calculation (cf.\
\cite{Clift}), we now come to the following conclusion.

\begin{prop}\ \ $N_R=0\ \Longleftrightarrow\ N_J=0$. \qed
\end{prop}

It is of some interest to recall here the following characterization of
$N_J=0$.

\begin{lemma}
$N_J=0$ if and only if for all basic vector fields $X,Y$, we have
\begin{equation}
\DH{JX}J(Y)-J(\DH{X}J(Y))=\DH{JY}J(X)-J(\DH{Y}J(X)). \label{NJzero}
\end{equation}
\end{lemma}

\proof\ \ Taking into account that $[X,Y]=\subH{[X,Y]}$ for basic
vector fields, we have
\begin{eqnarray*}
\lefteqn{[JX,JY]+J^2[X,Y]-J([JX,Y]+[X,JY]) = \DH{JX}(JY)-\DH{JY}(JX)} \\
&& \mbox{} +J^2(\DH{X}Y-\DH{Y}X) -J(\DH{JX}Y-\DH{Y}(JX)+\DH{X}(JY)-\DH{JY}(X))\\
&=& \DH{JX}J(Y)-J(\DH{X}J(Y))-\DH{JY}J(X)+J(\DH{Y}J(X)),
\end{eqnarray*}
which gives the desired result. \qed

To conclude this section, we look at an interesting case of
non-vanishing $\lie{\Gamma}R$. Observe first that $\Gamma$, by
construction, is the Hamiltonian vector field associated to the
energy function $E_L=\Delta(L)-L$, with respect to the symplectic
form $\omega_L=d\theta_L$, i.e.\ we have
\begin{equation}
i_\Gamma d\theta_L = - dE_L. \label{defgamma}
\end{equation}
A general theorem proved in \cite{CST} then implies that in the
present context, the fact that the 2-form $\omega_1$
(\ref{omega1}) is closed, is equivalent to stating that
\begin{equation}
i_{\lie{\Gamma}R}d\theta_L= -2\,dd_RE_L. \label{propCST}
\end{equation}
Since $R$ is symmetric with respect to $\omega_L$, the same is
true for $\lie{\Gamma}R$, so that invariance of $R$ is equivalent
to having $dd_RE_L=0$. The latter was the starting point for an
application of a bi-differential calculus in \cite{CST}, to which
we shall return in the next section. The more general situation of
a gauged bi-differential calculus in \cite{CST}, corresponds to
the assumption that for some basic function $f$,
\begin{equation}
dd_RE_L= df\wedge dE_L. \label{gauging}
\end{equation}
Via the equality (\ref{propCST}), this assumption is equivalent to
stating that (cf.\ Prop.~5.3 in \cite{CST})
\begin{equation}
\lie{\Gamma}R = \Gamma\otimes df - \xi_f\otimes dE_L,
\label{gauging2}
\end{equation}
where $\xi_f$ is the Hamiltonian vector field associated to $f$. It easily
follows, from $i_{S\xi_f}\omega_L=-S(i_{\xi_f}\omega_L)= S(df)=0$, that
$\xi_f$ is vertical, say $\xi_f=\V{X_f}$ for some $X_f$ along $\tau$. Then,
\[
i_{\xi_f}\omega_L (\H{Y}) = \K{g}(\V{X_f},\H{Y})= g(X_f,Y) = -
\H{Y}(f).
\]
Hence, in terms of fields along the projection $\tau$, $X_f$ is
defined by
\begin{equation}
g(X_f,Y) = - \H{Y}(f)= -\dh{}\!f(Y) \quad\mbox{or}\quad X_f\hook g =
-\dh{}\!f. \label{Xf}
\end{equation}

In the next section we specialize to the case that the Lagrangian
comes from a (pseudo) Riemannian metric on $Q$, and will focus
most of the attention on the characterization of so-called special
conformal Killing tensors in their tangent bundle manifestation.

\section{The Riemannian case}

Let $g$ be a symmetric, non-singular type $(0,2)$ tensor field on
$Q$ and put $L=\onehalf g_{ij}u^iu^j$. The non-linear connection
defined by the Euler-Lagrange equations then is the (linear)
Levi-Civita connection of $g$, i.e.\ the connection coefficients
(\ref{conn}) are of the form
\begin{equation}
\Gamma^i_j = \Gamma^i_{jk}u^k, \label{LCconn}
\end{equation}
where $\Gamma^i_{jk}$ are the classical Christoffel symbols.

It is important to understand first how the fundamental covariant
derivative operators of the calculus along $\tau$, as referred to
in the previous section, relate to classical tensor calculus in
this case. For example, if $J$ is a type $(1,1)$ tensor field on
$Q$, then both the dynamical covariant derivative $\del J$ and
horizontal covariant derivatives like $\DH{X}J$ relate to the
classical covariant derivative $J^i_{j|k}$, as follows:
\begin{equation}
(\del J)^i_j = J^i_{j|k}u^k, \qquad \left(\DH{X}J\right)^i_j =
J^i_{j|k}X^k. \label{covder}
\end{equation}
In fact, in this situation we have $\del= \DH{\T}$. It follows
that, in particular, $\del g=0$ and $\DH{X}g=0,\ \forall X$ (and
of course also $\DV{X}g=0$ because $g$ is basic). Another specific
property of the case of a quadratic spray is that the so-called
{\em deviation}, which in the language of the calculus along
$\tau$ is $\del \T$, is zero. It is actually of interest to list
all covariant derivatives of $\T$ here:
\begin{equation}
\del \T=0, \qquad \DV{X}\T=X, \qquad \DH{X}\T=0. \label{covT}
\end{equation}
These properties are easy to verify in coordinates, but let's take
the opportunity to mention also the general commutator identity
\begin{equation}
[\del,\DV{X}] = \DV{\del X} - \DH{X}, \label{deldv}
\end{equation}
which can be used to show in a coordinate free way that the first
two relations (\ref{covT}) imply the third.

We now look at the characterization of $R$ in this context, more particularly
the specification of the tensor field $U$. Note first that everything should
now be easily expressible in terms of the metric since we have, for example,
that
\begin{equation}
\theta_L(X)=g(\T,X)\quad\mbox{or}\quad \theta_L=\T\hook g. \label{newthetaL}
\end{equation}
We should keep in mind also that $L=E_L$ is a first integral, so that $\del
L=\Gamma(L)=0$, and that it further follows from (\ref{newthetaL}), using the
commutator property $[\del,\dv{}]=-\dh{}$, that
\begin{equation}
0=\del\theta_L=\del\dv{}\!L=-\dh{}\!L. \label{dhLzero}
\end{equation}
Now, concerning the determination of $U$, we have
\begin{eqnarray*}
\dh{}\!(J\theta_L)(X,Y) &=& \DH{X}(J\theta_L)(Y) - \DH{Y}(J\theta_L)(X) \\
&=&  \DH{X}(g(\T,JY)) - g(\T,J(\DH{X}Y)) - \DH{Y}(g(\T,JX)) +g(\T,J(\DH{Y}X)).
\end{eqnarray*}
Taking into account that $\DH{X}g=0$ and $\DH{X}\T=0$, we conclude that $U$ is
determined by
\begin{equation}
g(UX,Y) = g(\T,\DH{X}J(Y)-\DH{Y}J(X)) = g(\T,\dh{}\!J(X,Y)). \label{U}
\end{equation}

In coordinates, we of course work with the adapted frame of horizontal and
vertical vector fields on $TQ$ and their dual 1-forms, which are
\begin{equation}
\left\{H_i = \fpd{}{q^i} - \Gamma^k_i\fpd{}{u^k}, \ V_i =
\fpd{}{u^i}\right\}, \quad \left\{dq^i, \ \eta^j=du^j + \Gamma^j_k
dq^k\right\}. \label{adaptedframe}
\end{equation}
The coordinate expression for $R$ then becomes
\begin{equation}
R= J^i_j H_i\otimes dq^j + \ov{J}^i_j V_i\otimes\eta^j + U^i_j
V_i\otimes dq^j, \label{coordR}
\end{equation}
where
\begin{equation}
U^i_j = g^{ik}(J^m_{k|j}-J^m_{j|k})g_{ml}u^l. \label{coordU}
\end{equation}

It is interesting to have intrinsic expressions also, which implicitly
determine the vertical and horizontal covariant derivatives of $U$.

\begin{prop}
For $Z\in\vectorfields{\tau}$, $\DV{Z}U$ and $\DH{Z}U$ are
determined by
\begin{eqnarray}
g(\DV{Z}U(X),Y) &=& g(Z,\DH{X}J(Y)-\DH{Y}J(X)), \label{dvU} \\
g(\DH{Z}U(X),Y) &=& g(\T,\DH{}\DH{}J(Z,X,Y) - \DH{}\DH{}J(Z,Y,X)),
\label{dhU}
\end{eqnarray}
\end{prop}
\proof\ \ The proof is a straightforward computation, which starts from
(\ref{U}) and takes the properties (\ref{covT}) into account, remembering
further that
\[
\DH{}\DH{}J(Z,X,Y) = \DH{Z}\DH{X}J(Y) - \DH{\DH{Z}X}J(Y). \qquad
\qquad \mbox{\qed}
\]

Let us come back in this case of particular interest to the invariance of $R$,
or the more general assumption (\ref{gauging2}). If $h=\onehalf g^{ij}p_ip_j$
is the corresponding Hamiltonian, we know from the results of the preceding
section that $\lie{\Gamma}R=0$ is equivalent to $dd_{\tilde{J}}h=0$. It was
shown by a coordinate calculation in \cite{CST} that the latter condition is
equivalent to $J=\ov{J}$ and $\del J=0$. This is somewhat surprising now,
since Proposition~5, in general, imposes more conditions for having
$\lie{\Gamma}R=0$. Of course, if $\del J=0$, it follows from (\ref{covder})
that also $\DH{X}J=0$ and thus from (\ref{U}) that also $U=0$. This implies in
particular, from Proposition~2, that when $R$ is invariant, it must be equal
to $J^c$. But the more interesting information which follows from comparison
with Proposition~5 is that apparently, when $J$ is symmetric and parallel in
the Riemannian case, it will automatically commute with the Jacobi
endomorphism $\Phi$. This is not a trivial property to recognize. We therefore
propose to verify in Appendix~A by a direct calculation that it is indeed a
correct statement. That direct proof is of interest in its own right, because
it illustrates how one can proceed with an integrability analysis in this
context.

Now assume that (\ref{gauging2}) holds for some $f\in\cinfty{Q}$
and where $\xi_f=\V{X_f}$, with $X_f$ defined by (\ref{Xf}). In
the present situation, we further have $E_L=L$ and $\Gamma=\H{\T}$
(since $\del\T=0$).

\begin{thm} Under the present circumstances, the tensor field
$\lie{\Gamma}R$ is of the form (\ref{gauging2}) if and only if
$J=\ov J$ and further satisfies
\begin{equation}
\del J =\onehalf(\T\otimes \dh{}\!f - X_f\otimes \theta_L).
\label{scK}
\end{equation}
In addition, $U$ then is of the form
\begin{equation}
U= - \onehalf(\T\otimes \dh{}\!f + X_f\otimes \theta_L),
\label{scKU}
\end{equation}
and $J$ further has the property
\begin{equation}
\Phi J - J \Phi =  \onehalf(\T \otimes \del\dh{}\!f + \del
X_f\otimes\theta_L). \label{PhiJ}
\end{equation}
Finally, the tensor field $R$ itself then is given by
\begin{equation}
R=J^c - \Delta\otimes df. \label{scKR}
\end{equation}
\end{thm}
\proof\ \ The right-hand side of (\ref{gauging2}), when evaluated
on some $\V{X}$, results in $-\V{X}(E_L)\V{X_f} = -
\theta_L(X)\V{X_f}$. Comparison with (\ref{lgamR1}) shows that
this requires $J$ to be symmetric, plus the condition that
\[
U=-\del J - X_f\otimes \theta_L.
\]
Proceeding in the same way for an arbitrary horizontal argument
$\H{X}$, comparison with (\ref{lgamR2}) reveals the requirements
\[
U=\del J - \T\otimes \dh{}\!f,
\]
and
\[
\del U + \Phi J - J \Phi = - X_f\otimes \dh{}\!L=0.
\]
Compatibility of the two expressions for $U$ above, immediately leads to the
conclusions (\ref{scK}) and (\ref{scKU}). The first of these puts a
restriction on $J$, while the second in fact is then automatically satisfied.
To see this, we compute $\dv{}\!\del J$ from (\ref{scK}). Remember (see
\cite{MaCaSaI}) that for a vector valued 1-form along $\tau$, of the form
$\alpha\otimes X$, an exterior derivative such as $\dv{}$ is computed via the
rule: $\dv{}(\alpha\otimes X)= \dv{}\!\alpha \otimes X -
\alpha\wedge\dv{}\!X$. Now $\dv{}\!\dh{}\!f = - \dh{}\!\dv{}\!f=0$,
$\dv{}\!\T= I$ (the identity tensor), $\dv{}\!\dv{}\!L=0$ and finally also
$\dv{}\!X_f=0$ from (\ref{Xf}). It follows that
\begin{equation}
\dh{}\!J = \dv{}\!\del J = - \onehalf \dh{}\!f\wedge I. \label{dhJ}
\end{equation}
The defining relation (\ref{U}) for $U$ then easily leads to (\ref{scKU}) and
can in fact also be rewritten as
\begin{equation}
g(UX,Y) = - \onehalf\dh{}\!f\wedge\theta_L(X,Y). \label{Ubis}
\end{equation}
The final requirement that $\del U$ should be equal to $J\Phi-\Phi J$ leads
immediately, from (\ref{scKU}), to (\ref{PhiJ}), or can equivalently, from
(\ref{Ubis}), be expressed as
\begin{equation}
(\Phi J-J\Phi)\hook g = \onehalf \del\dh{}\!f\wedge \theta_L. \label{PhiJ2}
\end{equation}
The point is, however, that this again is not an extra condition, but a
consequence of the fundamental condition (\ref{scK}). To see this, recall that
(\ref{gauging2}) is equivalent to (\ref{gauging}), which in turn, when
translated into the corresponding cotangent bundle property reads $dd_{\tilde
J}h=df\wedge dh$. It was shown by a direct coordinate calculation in
\cite{CST} that this condition requires the symmetric $J$ to be a so-called
{\sl special conformal Killing tensor\/}. We shall verify in coordinates below
that this is exactly the condition (\ref{scK}). Hence, (\ref{PhiJ}) must be a
corollary and one could obtain it in a direct way by following the pattern of
the integrability analysis in Appendix~A.

The final statement (\ref{scKR}) about $R$ follows directly from comparison
between (\ref{Jcxv}-\ref{Jcxh}) and (\ref{Rvx}-\ref{Rhx}), knowing that $J=\ov
J$ and using (\ref{scK}) and (\ref{scKU}), with $\V{\T}=\Delta$. \qed

In coordinates, the condition (\ref{scK}) reads
\begin{equation}
J^i_{j|k} = \onehalf \left(\delta^i_k\fpd{f}{q^j} +
g^{il}\fpd{f}{q^l}g_{jk}\right). \label{scKcoord1}
\end{equation}
The more elegant coordinate expression is obtained by lowering an index and
reads
\begin{equation}
J_{lj|k} = \onehalf\left(g_{lk}\fpd{f}{q^j} + g_{jk}\fpd{f}{q^l}\right),
\label{scKcoord2}
\end{equation}
which is indeed the defining relation for a special conformal Killing tensor
as used in previous work (see e.g.\ \cite{CST,CraSar1}). An advantage of the
present framework is that we do obtain an easy to handle and elegant,
intrinsic expression also for the condition on the tensor $J$ in its type
$(1,1)$ appearance, which is after all the way in which $J$ is originally
conceived.

We finish this overview of the Riemannian case by briefly re-deriving the most
important properties of special conformal Killing tensors from (\ref{scK}).

\begin{thm} If $J$ is symmetric and satisfies (\ref{scK}) for some function $f\in\cinfty{Q}$,
then $N_J=0$ and $f={\rm tr\,} J$; moreover, if $J$ is non-singular, then its
cofactor tensor $A$ is a Killing tensor.
\end{thm}
\proof\ \ Acting with $\DV{X}$ on (\ref{scK}), and knowing that $\DV{X}X_f=0$,
it follows from the commutator property (\ref{deldv}) that
\begin{equation}
\DH{X}J = \onehalf (X\otimes\dh{}\!f - X_f\otimes \DV{X}\theta_L).
\label{dhscK}
\end{equation}
We then get
\begin{eqnarray*}
\lefteqn{\DH{JX}J(Y)-J(\DH{X}J(Y))-\DH{JY}J(X)+J(\DH{Y}J(X))} \\
&& = \onehalf \Big(\DV{JY}\theta_L(X)-\DV{JX}\theta_L(Y)\Big)X_f + \onehalf
\Big(\DV{X}\theta_L(Y)-\DV{Y}\theta_L(X)\Big)JX_f.
\end{eqnarray*}
The second term manifestly vanishes because
$\dv{}\!\theta_L=\dv{}\!\dv{}\!L=0$. Taking $X$ and $Y$ to be basic for
simplicity, the coefficient of $X_f$ can be re-written as
\[
\DV{JY}(\theta_L(X)) - \DV{JX}(\theta_L(Y))= \DV{JY}(g(\T,X)) -
\DV{JX}(g(\T,Y))= g(JY,X) - g(JX,Y),
\]
which is zero in view of the symmetry of $J$. Lemma~2 now implies $N_J=0$.

From (\ref{scK}), we get
\begin{eqnarray*}
\del({\rm tr\,}J) &=& \onehalf\Big(\langle\T,\dh{}\!f\rangle -
\langle X_f,\theta_L\rangle\Big) \\
&=& \onehalf(\del f - g(\T,X_f))= \del f.
\end{eqnarray*}
Hence, $f={\rm tr\,} J$ (up to a constant, which is irrelevant).

Finally, if $A$ is the cofactor of $J$, meaning that $JA=(\det
J)I$, we have
\begin{equation}
(\DH{X}J)A = - J(\DH{X}A) + \H{X}(\det J)\, I. \label{aux}
\end{equation}
Again, it suffices to let $X,Y,Z$ in what follows be basic vector fields, so
that for example $\DV{X}\theta_L(Z)=\DV{X}(\theta_L(Z))=g(X,Z)$. For the sake
of uniformity, we keep using the operators of the calculus along $\tau$,
although everything here of course happens on the base space $Q$ (and
expressions such as $\H{X}(f)$ mean simply $X(f)$). From (\ref{dhscK}), it
follows that
\[
\DH{X}J(AY) = \onehalf\Big(\H{(AY)}(f)X-g(X,AY)X_f\Big).
\]
Using this to compute $g(\DH{X}J(AY),AZ)$ and taking a cyclic sum
over $X,Y,Z$ (indicated by an ordinary summation symbol), we
readily obtain, knowing that also $A$ is symmetric,
\[
\sum g(\DH{X}J(AY),AZ) = \sum \H{(AX)}(f)\,g(Y,AZ).
\]
Next, using this to compute $\sum g(J(\DH{X}A)Y,AZ)= (\det J)\sum
g(\DH{X}A(Y),Z)$ via (\ref{aux}), we arrive at
\[
(\det J)\sum g(\DH{X}A(Y),Z)= \sum (\H{X}(\det J)- \H{(AX)}(f))
\,g(Y,AZ).
\]
But we know that $N_J=0$ implies that $d_J(\det J)=(\det J)\,d({\rm tr\,}J)$
(see e.g.\ \cite{BolMat,CraSar2}), which can be written as $\dh{}\!(\det
J)(JX^\prime) = (\det J)\,\dh{}\!f(X^\prime)$, for all $X^\prime$. Taking
$X^\prime=AX$, it follows that
\begin{equation}
\H{X}(\det J) = \H{(AX)}(f), \label{dJdetJ}
\end{equation}
which in turn leads to
\begin{equation}
\sum g(\DH{X}A(Y),Z) =0. \label{killing1}
\end{equation}
This is the way to express that $A$, as type $(1,1)$ tensor, is a
Killing tensor. The more familiar way is to look at the type
$(0,2)$ tensor field $\tilde{A}$, obtained by lowering an index,
so that $J\hook\tilde{A}=(\det J)\,g$. The cyclic sum condition
\begin{equation}
\sum \DH{X}\tilde{A}\,(Y,Z)=0 \label{killing2}
\end{equation}
then is equivalent to (\ref{killing1}). \qed

\section{An outlook for further study}

As stated in the introduction, our goal is to develop
generalizations of the classical cases of Hamilton-Jacobi
separable systems, or completely integrable systems, of which we
have not given any examples here, because such examples can
abundantly be found in the cited literature. There are reasons to
believe that a tangent bundle approach will then have advantages
over a cotangent bundle framework. The present study is a
preliminary investigation about understanding how things work on a
tangent bundle. Even so, we have already obtained in Sections~3
and 4 some general results relating to an arbitrary Lagrangian
function (not necessarily one of `mechanical type'). But for a
full generalization, also more general type $(1,1)$ tensors $J$
should be allowed, with components depending on coordinates and
velocities (or coordinates and momenta). For such a $J$, the
notion of complete lift is lost, so how to proceed? The point now
is that we indeed have an idea of how to proceed in the tangent
bundle set-up. It suffices to look at the expression (\ref{Jc})
for $J^c$ and to observe that the right-hand side is perfectly
defined also for a tensor field $J$ along the projection $\tau$.
In fact, we are then talking about a more general lifting
procedure, which has been fully developed already in
\cite{MaCaSaI,MaCaSaII} and is sometimes referred to as the
$\Gamma$-lift. For a given $J$ along $\tau$ which is not basic,
the formula
\begin{equation}
{\mathcal J}_\Gamma J = \H{J} + \V{(\del J)}, \label{gammalift}
\end{equation}
indeed defines a type $(1,1)$ tensor field on $TQ$, which depends on a given
second-order equation field $\Gamma$. A number of the calculations which will
be involved in such a generalization start off in exactly the same way as in
the present paper, but of course without the simplifications coming from
certain objects being basic.

A particular case of interest which could already significantly
generalize the well-known Riemannian situation of the previous
section, is to let $g$ be the metric along $\tau$ coming from a
Lagrangian which is the square of a Finsler function $F$. It is
then appropriate to start from a $J$ along $\tau$ (with the zero
section of $TQ$ excluded) which is homogeneous of degree zero in
the velocities. Thus, one can use the lifting procedure above,
where $\Gamma$ is the Euler-Lagrange field of $F^2$. Work along
these lines is in progress.

\section*{Appendix A: Aspects of integrability analysis}

Starting from an arbitrary $J$ on $Q$ and defining $R$ on $TQ$ by
(\ref{ourR}), the conditions for having $\lie{\Gamma}R=0$ are,
according to Proposition~5, that $J$ is symmetric and parallel,
and further commutes with the Jacobi endomorphism $\Phi$. But we
have argued indirectly in Section~5 that in the Riemannian case,
the third condition must be an automatic consequence of the first
two. An explicit verification of this fact can only come from an
integrability analysis on the partial differential equations
satisfied by $J$.

The assumption is that $\del J=0$, and since $J$ is basic, also
$\DV{X}J=0$ for all $X$, so that (\ref{deldv}) implies that also
$\DH{X}J$ will be zero for all $X$ (hardly a surprise of course in
view of (\ref{covder}). The next interesting commutator to look at
here is $[\del,\DH{X}]$, or more generally $[\DH{X},\DH{Y}]$. Its
general expression reads (see \cite{MaCaSaII})
\begin{equation}
[\DH{X},\DH{Y}]= \DH{\subH{[X,Y]}} + \DV{{\mathcal R}(X,Y)} +
\mu_{{\rm Rie}(X,Y)}. \label{dhdh}
\end{equation}
Clearly, only the last term matters here; it is a derivation
which, when acting on a vector field $Z$ along $\tau$ is given by
\[
\mu_{{\rm Rie}(X,Y)}Z= {\rm Rie}(X,Y)Z = - \DV{Z}{\mathcal
R}(X,Y).
\]
In fact ${\rm Rie}$ is here simply the classical Riemann tensor.
Recall also that we have the following relations linking $\Phi$
and ${\mathcal R}$,
\begin{equation}
\dv{}\!\Phi= 3{\mathcal R}, \qquad \Phi(X)={\mathcal R}(\T,X),
\label{PhiR}
\end{equation}
which implies for example that $\Phi(\T)=0$.

Now, to express that $[\del,\DH{X}]J$ must be zero, we compute
\[
[\del,\DH{X}](JY) - J([\del,\DH{X}]Y)= {\rm Rie}(\T,X)(JY) -
J({\rm Rie}(\T,X)Y).
\]
Using the second of (\ref{PhiR}) and the property $\DV{Y}\T=Y$, we
can write $\DV{Y}{\mathcal R}(\T,X)= \DV{Y}\Phi(X)+{\mathcal
R}(X,Y)$, by which the above expression in the end reduces to
\begin{equation}
A(X,Y):= - \DV{JY}\Phi(X)-{\mathcal R}(X,JY)+J(\DV{Y}\Phi(X)) +
J({\mathcal R}(X,Y)) =0. \label{aux2}
\end{equation}
In particular, knowing that $\Phi$ here is quadratic in the
velocities so that $\DV{\T}\Phi=2\Phi$, it follows that
\begin{equation}
A(X,\T):= - \DV{J\T}\Phi(X)-{\mathcal R}(X,J\T) + J(\Phi(X))=0,
\label{aux3}
\end{equation}
from which we further obtain that
\begin{equation}
0=g(A(X,\T),\T)=-g({\mathcal R}(X,J\T),\T) - g(\DV{J\T}\Phi(X),\T)
+ g(J\Phi(X),\T). \label{aux4}
\end{equation}
Now we recall from the study of the inverse problem of the
calculus of variations that also $\Phi$ is symmetric with respect
to $g$ (see e.g.\ Theorem~8.1 in \cite{MaCaSaII}). It then follows
from the Bianchi identity
\begin{equation}
\sum g({\mathcal R}(X,Y),Z)=0, \label{Bianchi}
\end{equation}
applied with arguments $X,JY,\T$, that
\begin{equation}
g({\mathcal R}(X,JY),\T)= g(\Phi J(Y),X) - g(\Phi X,JY) =0.
\label{aux5}
\end{equation}
So the first term on the right in (\ref{aux4}) vanishes. Moreover,
also $\DV{X}\Phi$ is symmetric with respect to $g$ for all $X$,
which implies that the second term can be re-written as $-
g(X,\DV{J\T}\Phi(\T))= g(X,\Phi(\DV{J\T}\T))=g(X,\Phi J(\T))$. The
conclusion from (\ref{aux4}) therefore is that $g(J\Phi(X),\T)=0$.
Computing the $\DV{Y}$-derivative of this result and subtracting
the same expression with $X$ and $Y$ interchanged, it follows by
using the first property in (\ref{PhiR}) and the symmetry of $J$
and $\Phi$ that
\begin{equation}
3\,g({\mathcal R}(X,Y),J\T) = g((J\Phi-\Phi J)X,Y). \label{aux6}
\end{equation}
Next, we use the Bianchi identity again to write
\[
g({\mathcal R}(X,Y),J\T) = g({\mathcal R}(X,J\T),Y) - g({\mathcal
R}(Y,J\T),X),
\]
and make use of (\ref{aux3}) to compute the right-hand side.
Taking $X$ and $Y$ to be basic for simplicity, we can write
$g(\DV{J\T}\Phi(X),Y)=\DV{J\T}(g(\Phi X,Y))$ and likewise for the
term with $X$ and $Y$ interchanged. It then readily follows that
$g({\mathcal R}(X,Y),J\T) = g((J\Phi-\Phi J)X,Y)$. Comparison with
(\ref{aux6}) leads to the conclusion that both sides must be zero,
for arbitrary $X,Y$. Hence, we have shown in a direct way that in
the Riemannian case,
\begin{equation}
J=\ov{J} \quad \mbox{and}\quad \del J=0 \quad \Rightarrow\quad
\Phi J=J\Phi. \label{[PhiJ]}
\end{equation}
We had an indirect proof of this fact in Section~5. It is of
interest to illustrate that it is indeed a non-trivial property by
looking at coordinate expressions.

We have
\begin{equation}
\Phi^i_j = {\mathcal R}^i_{kj}u^k,\quad {\mathcal
R}^i_{kj}=R^i_{ljk}u^l, \quad \mbox{and thus} \quad \Phi^i_j=
R^i_{ljk}u^ku^l, \label{coordPhi}
\end{equation}
where $R^i_{ljk}$ are the components of the Riemann tensor, and
are skew-symmetric in the last two subscripts. Now, from
$J^i_{j|k}=0$, taking a further covariant derivative, swapping
indices and using the Ricci identities, the property which
immediately follows is $J^i_j R^j_{kml}=R^i_{jml}J^j_k$. But the
commutation of $J$ and $\Phi$ is a different property and says
that $J^i_j(R^j_{kml}+R^j_{lmk})=(R^i_{kjl}+R^i_{ljk})J^j_m$.

\section*{Appendix B: Darboux-Nijenhuis coordinates}

It is well known that on a general (regular) Poisson-Nijenhuis manifold of
dimension $2n$, if the recursion operator $R$ has $n$ distinct eigenvalues,
their exist so-called Darboux-Nijenhuis coordinates, which diagonalize $R$ and
are at the same time Darboux coordinates for the symplectic form (see e.g.\
\cite{Turiel,FP}). This will apply in particular to the general situation on
$TQ$, described in Sections~3 and 4. We wish to investigate here in some
detail what the structure is of the eigenspaces of our $R$ and how the
construction of Darboux-Nijenhuis coordinates works when the eigenvalues are
maximally distinct.

We begin by establishing results which are valid without special assumptions
on the type $(1,1)$ tensor $J$ on $Q$, except that we will only consider real
eigenvalues.

\begin{lemma} If $\xi=\H{X}+\V{Y}$ is an eigenvector of $R$,
corresponding to the eigenvalue $\lambda$, then
\begin{equation}
JX=\lambda X \quad\mbox{and}\quad UX+\ov{J}Y=\lambda Y.
\label{eigenvector}
\end{equation}
It follows in particular that $X$ is an eigenvector of $J$.
\end{lemma}
\proof\ \ Using the characterization of $R$ as described by
(\ref{Rvx}-\ref{Rhx}), it is immediate to see that $R\xi=\lambda \xi$ is
equivalent to the two relations (\ref{eigenvector}). \qed

\begin{lemma} $J$ and $\ov{J}$ have the same eigenvalues. In fact,
if $X$ is an eigenvector of $J$, then $X\hook g$ is an eigenform
of $\ov{J}$ with the same eigenvalue.
\end{lemma}
\proof\ \ We have $\ov{J}^l_j=g^{lk}J^i_k g_{ij}$, and therefore
\[
\ov{J}^l_j - \lambda \delta^l_j = g^{lk}(J^i_k-\lambda
\delta^i_k)g_{ij}.
\]
Both statements now easily follow. \qed

\begin{lemma} Suppose that $J$ is non-degenerate and has $n$
distinct eigenvalues (which then are non-zero). Then, if
$JX=\lambda X$, there exists a vector field $Y$ along $\tau$, such
that $\ov{J}Y=\lambda Y - UX$.
\end{lemma}
\proof\ \ From $g(JX,Y)=g(X,\ov{J}Y)=\lambda g(X,Y)$, it follows
that $g(X,\ov{J}Y-\lambda Y)=0, \forall Y$. Extending $X=X_1$ to
an orthogonal frame $\{X_1,\ldots,X_n\}=\{X_1,X_\alpha\}$ for $g$,
and putting $Y=a^iX_i$, it follows that $a^i(\ov{J}X_i-\lambda
X_i)\in {\rm sp\,}\{X_\alpha\}, \forall a^i$, which implies that
\[
\ov{J}X_i= \lambda X_i + b^\alpha_iX_\alpha, \qquad i=1,\ldots,n
\]
for some functions $b^\alpha_i$. We know that $\lambda$ is an eigenvalue of
$\ov{J}$ as well, and that its eigenvalues are distinct. Hence, there exists a
unique vector field of the form $X_1+c^\alpha X_\alpha$ which spans the kernel
of $\ov{J}-\lambda I$. But $(\ov{J}-\lambda I)(X_1+c^\alpha
X_\alpha)=(b^\alpha_1+c^\beta b^\alpha_\beta)X_\alpha$, so the fact that
unique functions $c^\beta$ exist which make this zero implies that
$\det(b^\alpha_\beta)\neq 0$. Now consider the equation $\ov{J}Y=\lambda Y -
UX$ for the unknown $Y=a^iX_i\in\vectorfields{\tau}$. Since $g(UX,Y)$ is
skew-symmetric in $X,Y$, we know that $g(UX_1,X_1)=0$ and thus $UX_1=d^\alpha
X_\alpha$ for some functions $d^\alpha$. The equation for $Y$ can now be
written in the form $a^\beta b^\alpha_\beta= - d^\alpha - a^1 b^\alpha_1$ and
clearly has a unique solution for the $a^\beta$ for each arbitrary choice of
$a^1$. \qed

\begin{prop} Let $J$ be diagonalizable with distinct non-zero eigenvalues.
Then a complete set of eigenvectors of $R$ can be constructed as follows: (i)
let $X_i$ denote the eigenvector of $J$ with eigenvalue $\lambda_i$ and $Z_i$
the eigenvector of $\ov{J}$ with the same eigenvalue; (ii) for each $X_i$,
construct a vector $Y_i$ such that $\ov{J}Y_i=\lambda_iY_i - UX_i$. Then
$\V{Z_i}$ and $\H{X_i}+\V{Y_i}$ are eigenvectors of $R$, corresponding to the
eigenvalue $\lambda_i$.
\end{prop}
\proof\ \ We have
\begin{eqnarray*}
R(\V{Z_i}) &=& \V{(\ov{J}Z_i)} = \lambda_i\V{Z_i}, \\
R(\H{X_i}+\V{Y_i}) &=& \H{(JX_i)} + \V{(UX_i)}+\V{(\ov{J}Y_i)} =
\lambda_i(\H{X_i}+\V{Y_i}),
\end{eqnarray*}
from which the result follows. \qed

This much for the purely algebraic aspects. Now let us further
assume that $N_J=0$. Darboux-Nijenhuis coordinates in fact should
do three things at the same time: not only diagonalize $R$ in
coordinates, but also separate it, and bring the symplectic form
into canonical form. It was proved in the fundamental paper of
Fr\"olicher and Nijenhuis \cite{FrNi} that if $J$ is
(algebraically) diagonalizable and the eigenvalues have constant
multiplicity, then the {\sl necessary and sufficient condition for
diagonalizability in coordinates\/} is that ${\mathcal H}_J=0$,
where the {\sl Haantjes tensor\/} ${\mathcal H}_J$ can be defined
by
\begin{equation}
{\mathcal H}_J(X,Y)= J^2N_J(X,Y) + N_J(JX,JY) -
JN_J(JX,Y)-JN_J(X,JY). \label{haantjes}
\end{equation}
Obviously, $N_J=0$ implies ${\mathcal H}_J=0$, but evaluating
$N_J$ on eigenvectors $X$ and $Y$ belonging to different
eigenvalues, $\lambda, \mu$ say, further gives
\[
0=N_J(X,Y)= (\lambda -\mu)(X(\mu)Y+Y(\lambda)X),
\]
so that $X(\mu)=Y(\lambda)=0$. Hence, in coordinates which diagonalize $J$,
the eigenvalues will only depend on the coordinates of the corresponding
eigendistribution, which is the meaning of saying that $J$ is separable in
coordinates. Conversely, if $J$ is separable, one can verify in such
coordinates that $N_J=0$. In other words, $N_J=0$ (for a $J$ which has the
algebraic properties stated above) is the {\sl necessary and sufficient
condition for separability in coordinates\/}. Note in passing that the tools
for studying such issues when $J$ would more generally be a tensor field along
$\tau$ have been developed in \cite{MaCaSaIII}.

To understand what happens with $R$ on $TQ$ now, we need to look
at the expression of $R$ in a coordinate basis, rather than in the
adapted frame as in (\ref{coordR}); it reads (still for the
general situation described by Proposition~1)
\begin{equation}
R= J^i_j\fpd{}{q^i}\otimes dq^j + \ov{J}^i_j\fpd{}{u^i}\otimes
du^j + (U^i_j+\ov{J}^i_k\Gamma^k_j -
J^k_j\Gamma^i_k)\fpd{}{u^i}\otimes dq^j. \label{Rcoord1}
\end{equation}
The following procedure now will lead to Darboux-Nijenhuis
coordinates. First perform the Legendre transform $(q,u)
\rightarrow (q,p=\partial L/\partial u)$. Even though this is to
be regarded here as a change of coordinates on $TQ$, the result
will be that $R$ acquires the form of the complete lift
$\tilde{J}$ on $T^*Q$, i.e.\
\begin{equation}
R=J^i_j\left(\fpd{}{q^i}\otimes dq^j + \fpd{}{p_j}\otimes
dp_i\right) + p_k\left(\fpd{J^k_i}{q^j}-\fpd{J^k_j}{q^i}\right)
\fpd{}{p_i}\otimes dq^j. \label{Rcoord2}
\end{equation}
The 2-form $\omega_L$ meanwhile will already take its canonical
form in the $(x,p)$ coordinates. Now, assuming that $J$ has
distinct eigenvalues and zero Nijenhuis torsion, we know that
there exists a coordinate change on $Q$ which will diagonalize the
expression $J^i_j(\partial/\partial q^i \otimes dq^j)$ in such a
way that the eigenvalues depend on at most one new coordinate
$q'$. The resulting point transformation $(q,u)\rightarrow
(q',u')$ on $TQ$, when expressed in the non-tangent bundle
variables $(q,p)$, formally is a `canonical transformation'
$(q,p)\rightarrow(q',p')$, i.e.\ it defines another Darboux chart
for the symplectic form $\omega_L$ and it will have the additional
effect of diagonalizing (and separating) $R$.

From a tangent bundle point of view, the first step in this procedure is
rather unnatural, because it is not a tangent bundle change of coordinates. At
first sight, it may look like one should nevertheless not change the order of
the operations, because even though $J$ and $\ov{J}$ have the same
eigenvalues, a coordinate transformation which diagonalizes $J$ will generally
not at the same time diagonalize $\ov{J}$. However, the two coordinate changes
under consideration here are of course of a quite special type: a Legendre
transformation which does not change $q$ but changes the fibre coordinate, and
a point transformation. It is clear that such coordinate changes commute, so
one can just as well diagonalize $J$ first and then the subsequent Legendre
transform will not destroy the diagonal form of $J$, will bring $\omega_L$ in
canonical form, and at the same time will take care of the diagonalization of
$\ov{J}$.

That the reversed procedure is somewhat more natural for the
tangent bundle set-up may become clear in the special case that
$J$ is symmetric. It then follows from $g_{ij}J^i_k=g_{ik}J^i_j$
that in coordinates which diagonalize $J$, we will have
$g_{kj}(\lambda^{(k})-\lambda^{(j)})=0$ and thus $g_{kj}=0$ for
$j\neq k$. This gives useful information also when there is no
urge to pass to Darboux-Nijenhuis coordinates: it means that in
coordinates which diagonalize $J$, the given Lagrangian will
separate with respect to the velocity variables, i.e.\ it will
become of the form: $L=\sum_i L^i(q,u^i)$, where $L^i$ depends on
$u^i$ only.

\end{document}